\documentclass{amsart}
\usepackage{amsthm,amsmath,amsfonts,amssymb}
\usepackage{verbatim}
\usepackage{upgreek}
\usepackage{amsthm,amsmath,amsfonts}
\usepackage{mathrsfs}
\usepackage[colorlinks=true,linktocpage=true]{hyperref}
\usepackage{tocvsec2}
\usepackage{accents}
\usepackage{hyperref}
\usepackage{upgreek}
\usepackage{stmaryrd}
\usepackage[usenames,dvipsnames]{xcolor}
\hypersetup{urlcolor=cyan, citecolor=cyan, linkcolor=red}

\usepackage{framed}
\usepackage{bbm,dsfont}
\usepackage{color}
\usepackage{hyperref}
\numberwithin{equation}{section}
\numberwithin{table}{section}
\font\tenscrpt=eusm10
\font\sevenscrpt=eusm10 scaled 700
\font\fivescrpt=eusm10 scaled 500
\newfam\eusmfam
\textfont\eusmfam=\tenscrpt
\scriptfont\eusmfam=\sevenscrpt
\scriptscriptfont\eusmfam=\fivescrpt
\def\scr#1{{\fam\eusmfam\relax#1}}

\newtheorem{thm}{Theorem}[section]
\newtheorem{cor}{Corollary}[section]
\newtheorem{lem}{Lemma}[section]

\newtheorem{prop}{Proposition}[section]

\theoremstyle{definition}
\newtheorem{defn}{Definition}[section]

\newtheorem{rem}{Remark}[section]
\newtheorem{notn}{Notation}[section]
\newcommand{\thmref}[1]{Theorem~\ref{#1}}

\newcommand{\lemref}[1]{Lemma~\ref{#1}}

\newcommand{\prpref}[1]{Proposition~\ref{#1}}
\newcommand{\defnref}[1]{Definition~\ref{#1}}

\def\qed{\quad\vcenter{\hrule\hbox{\vrule height.6em\kern.6em\vrule}\hrule}}
\newenvironment{pf}{{\bigskip\textit{\newline Proof.}\quad}}{$\qed$\bigskip\newline}
\newenvironment{pf*}[1]{{\bigskip\textit{\newline#1.}\quad}}{$\qed$\bigskip\newline}
\def\ds{\displaystyle}

\def\lftdt{\left.}

\def\psxy{{K}^{\text{\tiny{\sc{BM}}}^d}_{s;x,y}}

\def\ptzs{{K}^{\text{\tiny{\sc{BM}}}}_{t;0,s}}

\def\K{{\mathbb K}}

\def\KBtxy{{\K}^{\text{\tiny{\sc{BTBM}}}^d}_{t;x,y}}
\def\KBtyo{{\K}^{\text{\tiny{\sc{BTBM}}}}_{t;0,y}}
\def\KBtzo{{\K}^{\text{\tiny{\sc{BTBM}}}}_{t;0,z}}

\def\KBtzmut{{\K}^{\text{\tiny{\sc{BTBM}}}}_{t;\mu t,z}}

\def\Xx{X^{x}}
\def\XxB{{\mathbb X}^{x}_{B}}

\def\XB{{\mathbb X}_{B}}
\def\XBt{{\mathbb X}_{B}(t)}
\def\XxBt{{\mathbb X}_{B}^{x}(t)}

\def\XxBz{{\mathbb X}_{B}^{x}(0)}
\def\EBTP{{\mathbb X}^{x,\infty}_{B,e}(t)}

\def\kEBTP{{\mathbb X}^{x,k}_{B,e}(t)}
\def\RBM{|B(t)|}
\def\XtBt{\tilde{\mathbb X}_{B}(t)}

\def\P{{\mathbb P}}
\def\Pt{\tilde{\mathbb P}}
\def\EP{{\mathbb E}_{\P}}

\def\E{{\mathbb E}}
\def\N{{\mathbb N}}

\def\Rd{{\mathbb R}^{d}}

\def\R{\mathbb R}

\def\Bxt{{\mathbb B}^{x}_{t}}

\def\X{\mathbb X}
\def\Rp{{\R}_+}

\def\sF{{\mathscr F}}
\def\sG{{\mathscr G}}
\def\sFt{{\mathscr F}_t}

\def\sGtB{{\mathscr G}_t^{B}}

\def\OFP{(\Omega,\sF,\P)}

\def\OFtPt{({\Omega},{\sF}_{t},{\Pt_{t}})}

\def\XitmuXB{\Xi_{t}^{\mu,\X_{B}}}
\def\XitmuX{\Xi_{t}^{\mu,B}}

\def\pa{\partial}

\def\df#1#2{\ds{\frac{#1}{#2}}}

\def\lbl#1{\label{#1}}

\def\pa{\partial}

\def\lang{\left<}
\def\rang{\right>}
\def\lab{\left|}
\def\rab{\right|}
\def\lpa{\left(}
\def\rpa{\right)}
\def\lbk{\left[}
\def\rbk{\right]}
\def\lbr{\left\{}
\def\rbr{\right\}}
\def\bdf{\begin{defn}}
\def\edf{\end{defn}}
\def\bcr{\begin{cor}}
\def\ecr{\end{cor}}
\def\bprop{\begin{prop}}
\def\eprop{\end{prop}}
\def\bnt{\begin{notn}}
\def\ent{\end{notn}}
\def\brm{\begin{rem}}
\def\erm{\end{rem}}
\def\blm{\begin{lem}}
\def\elm{\end{lem}}
\def\bpf{\begin{pf}}
\def\bpfs{\begin{pf*}}
\def\epf{\end{pf}}
\def\epfs{\end{pf*}}
\def\bfr{\begin{framed}}
\def\efr{\end{framed}}
\def\bitm{\begin{itemize}}
\def\eitm{\end{itemize}}
\def\ben{\begin{enumerate}}

\def\rencomrom{\renewcommand{\labelenumi}{(\roman{enumi})}}
\def\een{\end{enumerate}}

\def\scr{\mathscr}
\def\eqrf{\eqref}
\def\beq{\begin{equation}}
\def\beqs{\begin{equation*}}
\def\eeq{\end{equation}}
\def\eeqs{\end{equation*}}
\def\bsp{\begin{split}}
\def\esp{\end{split}}
\def\bc{\begin{cases}}
\def\ec{\end{cases}}
\def\bt{\begin{tabular}}
\def\et{\end{tabular}}

\def\bthm{\begin{thm}}
\def\ethm{\end{thm}}
\def\bpr{\begin{prop}}
\def\epr{\end{prop}}
\def\babs{\begin{abstract}}
\def\eabs{\end{abstract}}

\def\sCBxd{\scr{C}^{x,d}_{\mathrm{BTBM}}}

\def\ig{\iffalse}
\def\Var{\mbox{Var}}
\def\Cov{\mbox{Cov}}
\def\sgn{\mbox{sign}}
\def\ind{\mathbbm{1}}

\def\l{\ell}
\def\sB{\mathscr B}
\def\dira{{\mathfrak d}_{a}}
\def\dirz{{\mathfrak d}_{0}}
\def\eqdis{\stackrel {\mbox{\tiny{d}}}{=}}
\title[Brownian-time Change of Measure]{Brownian-time Change of measure}
\author{Bobomurod Abdurakhmanov and Hassan Allouba}
\address{Department of Mathematical Sciences, Kent State University, Kent, Ohio 44242
\\email: allouba@math.kent.edu}
\thanks{}
\subjclass[2000]{60H20, 60H15, 45H05, 45R05, 35R60, 60H30, 60J45, 60J35, 60J60, 60J65.}
\keywords{Brownian-time processes, change of measure.}

\begin{document}
\begin{abstract}
We prove a fundamental change of measure theorem for the Brownian-time Brownian motion and its associated Brownian-time processes class introduced by Allouba and Zheng in 2001.  This result, together with Allouba's prior work on (1) Brownian-time processes and their PDEs/SPDEs links and on (2) change of measure for SPDEs, is a critical building block in analyzing the behaviors of SDEs and SPDEs---of different types and orders---driven by Brownian-time noises and their relatives.
\end{abstract}
\maketitle
\tableofcontents
\section{Introduction and statement of main result}
\subsection{Genesis and motivation}
This is the first of several articles in the works, and it is based on the program started by Allouba in an unpublished work  \cite{Abtcom0}.   That work \cite{Abtcom0}, in turns, lies at the intersection of two major components of Allouba's research: 
\ben
\item
 his work on the Brownian-time processes and their connection to higher order and fractional PDEs/SPDEs \cite{AT}--\cite{Abtbs} and \cite{Abtpspde}--\cite{AZ01}, and 
 \item his work on space-time change of measure theorems and their application to different type and different order SPDEs \cite{AX16,Alksspde,Acom2,Acom1,Acom}.
 \een    

The concept of change of measure is an important ingredient and a powerful tool in stochastic analysis and its applications.   The fundamental idea and its applications to Gaussian and martingale processes may be traced back to Cameron and Martin \cite{CM44} and Girsanov \cite{Gir60}.   These applications range from the theoretical investigation of solutions to SDEs and, more recently,  SPDEs (existence, uniqueness, law smoothness, etc\ldots, see for example Allouba's work \cite{AX16,Alksspde,Acom2,Acom1,Acom}) to the famous arbitrage analysis in mathematical finance theory, among others.   

In this note, we give a fundamental change of measure theorem in the Brownian-time processes setting, which we now recall.
In \cite{AZ01,Abtp2}, Allouba and Zheng introduced a large and rich class they called Brownian-time processes, the simplest of which is simple Brownian-time Brownian motion (simple BTBM): $$\XxB:=\lbr\XxBt=\Xx\lpa\lab B_{t}\rab\rpa,t\ge0\rbr,$$ where $\Xx$ is a Brownian motion starting at $x\in\Rd$ and $B$ is an independent one dimensional Brownian motion starting at $0$ on a probability space $\OFP$.  In Allouba and Zheng \cite{AZ01}, they introduced a BTBM-based class that includes important processes such as Le Gall's Brownian snake and Burdzy's
iterated Brownian motion (IBM) and introduce many new ones. This class of Brownian-time processes generalizes the simple BTBM above in several directions to obtain: 
\begin{enumerate}
\renewcommand{\labelenumi}{(\alph{enumi})}
\item the Brownian snake\ig\footnote{The Brownian snake is the only Markov process in $\sCBxd$.}\fi---when $|B(t)|$ increases we generate a new independent path.  See \cite{LeGall3} for applications to the nonlinear PDE $\Delta u=u^2$;
\item $k$-Excursions-based Brownian-time Brownian motions ($k$-EBTBMs).  Here, we let $X^{x,1}, \ldots,X^{x,k}$ be independent copies of $\Xx$ starting from point $x$.  On each excursion interval of $\RBM$ use one of the
$k$ copies chosen at random.  We denote such a process by $\kEBTP$.  Of course, when $k=1$ we obtain a simple BTBM, and
when $k=2$ we obtain the iterated Brownian-motion (IBM) of Burdzy \cite{BurdzyVar};
\item $\infty$-Excursions-based Brownian-time Brownian motions ($\infty$-EBTBMs).  Use an independent copy of $X^x$ on each excursion interval of $\RBM$.
This is the $k\to\infty$ weak limit of (b) (a proof of this is given in \cite{AZ01}).  It is intermediate between the processes in (b) and the Brownian snake in (a).
Here, we go forward on a new independent path only after $|B(t)|$
reaches $0$.  This process is abbreviated as $\infty$-EBTP and is denoted by $\EBTP$.
\end{enumerate} 
Although the elements in this Brownian-time processes  contains many different stochastic processes, they all share the same one dimensional distribution; i.e., they all have the same density $\P\lbk\XxBt\in dy\rbk=\KBtxy dy$ where $\KBtxy$ is the BTBM density\footnote{For more on the regularizing properties of the BTBM density and the related linear Kuramoto-Sivashinsky kernel and their implications to PDEs, SPDEs, and stochastic integral equations, we refer the reader to \cite{Abtpspde}--\cite{AZ01}.}:
\beq\lbl{BTBMden}  
\KBtxy=2\int_0^\infty \psxy\ptzs ds=2\int_0^\infty\df{e^{-|x-y|^2/2s}}{(2\pi s)^{{d}/{2}}}\df{e^{-s^2/2t}}{\sqrt{2\pi t}}ds.
\eeq 
\ig
\bdf[The class $\sCBxd$ and the BTBM random variable]  We say that a stochastic process $\XxB$ is a Brownian-time Brownian motion starting at $x\in\Rd$, and write $\XxB\in\sCBxd$, iff its one dimensional distribution $\P\lbk\XxBt\in dy\rbk=\KBtxy dy$ for all $t>0$ and $\P\lbk\XxBz=x\rbk=1$ for all $x\in\Rd$.  If $x\in\Rd$, we call a $d$-dimensional random vector $\Bxt=\lpa{\Bxt}^{(1)},\ldots,{\Bxt}^{(d)}\rpa$---where  each ${\Bxt}^{(i)}$ is a one dimensional BTBM random variable with parameters $x^{(i)}$ and $t$ and the components ${\Bxt}^{(i)}$'s are all independent---a BTBM with parameters $x$ and $t$ and write $\Bxt\stackrel {\mbox{\tiny{d}}}{=}\mbox{BTBM}^{d}(x,t)$ iff $\KBtxy$ is the $\Bxt$ probability density\footnote{By $\stackrel {\mbox{\tiny{d}}}{=}$ we mean equal in distribution.}.  
\edf  
\fi
\subsection{The main change of measure result}
We state the main result in terms of BTBM-related quantities $X$ and $B$, and its quartic variation.  First, we need the definition of $p$-th variation.
\bdf[$p$-th variation]\lbl{def:pvar} 
Fix $t>0$ and let $\Pi$ be a partition of $[0,t]$, $\Pi=\{0=t_{0}\le  t_{1}\le\cdots\le t_{n}=t\}$, with mesh $|\Pi|=\max_{1\le k\le n}|t_{k}-t_{k-1}|$.  Then for any stochastic process $Z$, we call 
\beq\lbl{eq;pthvoverPi}
V_{t}^{(p)}(Z,\Pi):=\sum_{k=1}^{n}\lab Z(t_{k})-Z(t_{k-1})\rab^{p}
\eeq
the $p$-th variation of $Z$ over the partition $\Pi$.  The $p$-th variation of $Z$, $\lang Z\rang^{(p)}_t$, $p>0$, is then given by 
\beq
\lang Z \rang^{(p)}_t:=\lim_{|\Pi|\searrow0}\sum_{k=1}^{n}\lab Z(t_{k})-Z(t_{k-1})\rab^{p}=\lim_{|\Pi|\searrow0}V_{t}^{(p)}(Z,\Pi),
\eeq
where the limit is in probability\footnote{As remarked by Revuz and Yor \cite{RevuzYor}, we don't take the sup over all partition as in the analytic definition.  Also, one may also require the limit to exist in other senses as well (almost sure, $L^{q}$ for $q>0$, etc\ldots).}.
\edf  
Our main result gives a probability measure $\Pt_{t}$, for each $t>0$, under which a translated Brownian-time Brownian motion  with drift has a zero mean BTBM density\footnote{For simplicity of notation and exposition, we state our result for the one dimensional BTBM case but it holds for the $d$-dimensional case without any major changes to our proofs.}.  
\bfr
\bthm[Brownian-time change of measure]\lbl{thm:BTCOM}  Let $\XB=\lbr\XB(t),\sFt;0\le t<\infty\rbr$ be a $1$-dimensional BTBM starting at $0\in\R$ on the probability space $\OFP$, and let
\beq\lbl{btbmwdrft}
\XtBt:=\XBt-\mu t; \quad \mu\in\R,\ t\ge0,
\eeq 
be the BTBM $\XB$ shifted by $\mu t$.  Define $\Pt_{t}$ by the Radon-Nikodym derivative
\beq\lbl{eq:RN}
\df{d\Pt_{t}}{d\P}:=\XitmuXB=\exp\lbr\df{t}{\RBM}\lpa\mu\XB(t)-\frac16\mu^{2}\lang\XB\rang^{(4)}_t\rpa\rbr.
\eeq 
Then, for each fixed $t$, $\Pt_{t}$ is a probability measure on $(\Omega,\sFt)$ and $\XtBt$ has a mean-zero BTBM density at $t$ on $\OFtPt$$:$ $\Pt_{t}[\XtBt\in dz]=\KBtzo dz$.
\ethm
\efr   
\brm
\ben\rencomrom
\item Since the quartic variation of a BTBM is $\lang\XB\rang^{(4)}_t=3t$ (see \prpref{prp:BTBMqv} below), the Radon-Nikodym derivative may be rewritten as
\beq\lbl{eq:RNexpl0}
\XitmuXB=\exp\lbr\df{t}{\RBM}\lpa\mu\XB(t)-\frac12\mu^{2}t\rpa\rbr.
\eeq
\item The importance of this result is that it is fundamental in the use of BTBM change measure to analyze many SDEs and SPDEs that are driven by BTBM noise in time \cite{Aks,Abtp2,AZ01} and driven in space-time by Allouba'a Brownian-time Brownian sheet \cite{Abtbs}.  This, in turns, is exciting since BTBM is connected to a large class of higher order and fractional PDEs and SPDEs \cite{Allouba2023}--\cite{Abtbs} and \cite{Abtpspde}--\cite{AZ01}.  Moreover, BTBMs do not fall into any of the classical theory of stochastic processes since they are neither Gaussian, nor semimartingale, nor Markov processes.  
\item Since the marginal density is the same for all members of the Allouba-Zheng class of Brownian-time processes above, the proof for the above change of measure theorem for the BTBM class above is a minor notational adaptation of our proof presented here for the simple BTBM.  
\item Moreover, the explicit Radon-Nikodym derivative spells out explicitly the role of the Brownian clock, $|B_{t}|$ in the BTBM change of measure theory.
\item In contrast to the BTBM case, the corresponding Radon-Nikodym in the Brownian motion $B$ case is simply
\beq\lbl{eq:BMcase}
\XitmuX=\exp\lbr\lpa\mu B(t)-\frac12\mu^{2}t\rpa\rbr.
\eeq
Comparing \eqref{eq:RNexpl0} with \eqref{eq:BMcase}, we see that the Brownian-time of the BTBM manifests itself as a time-normalizing factor in the Radon-Nikodym derivative for the BTBM case. 
\item In the $d$-dimensional case, it is trivial to extend the Radon-Nikodym derivative in \eqref{eq:RN} to
\beq\lbl{eq:RNd}
\df{d\Pt_{t}}{d\P}:=\exp\lbr\df{t}{\RBM}\lpa\sum_{i=1}^{d}\mu^{(i)}\XB^{(i)}(t)-\sum_{i=1}^{d}\lpa\mu^{(i)}\rpa^{2}t\rpa\rbr
\eeq 
\een
\erm
\section{Proof of Result}
The key idea---throughout our proofs below---is that of conditioning on the inner Brownian motion $B$ in the BTBM and use its independence from the outer Brownian motion $X$, as in Allouba-Zheng and Allouba \cite{AZ01,Abtp2}.  We then use suitable stochastic analytic and Gaussian estimates to complete these proofs.

The proof of \thmref{thm:BTCOM} follows from the next two propositions.  
\ig
\bprp[$\Pt_{t}$ is a probability measure]\lbl{lem:a probab measure} The measure $\Pt_{t}$ given by 
\beq\lbl{eq:RNexpl}
\df{d\Pt_{t}}{d\P}:=\XitmuXB_{t}=\exp\lbr\df{t}{\RBM}\lpa\mu\XB(t)-\frac12\mu^{2}t\rpa\rbr
\eeq 
is a probability measure on $(\Omega,\sFt)$
\eprp
\bpf
\epf
\fi
\subsection{The main result sans the BTBM quartic variation}\lbl{sec:stripQV}
The following proposition is a restatement of the main result \thmref{thm:BTCOM}, expressed without the BTBM quartic variation.
\bprop[Under $\Pt_{t}$, $\XtBt$ is a zero-mean BTBM at each fixed $t$]\lbl{lem:crux}
Let $\XtBt$ be as given in \thmref{thm:BTCOM} and let
\beq\lbl{eq:RNexpl}
\df{d\Pt_{t}}{d\P}:=\XitmuXB=\exp\lbr\df{t}{\RBM}\lpa\mu\XB(t)-\frac12\mu^{2}t\rpa\rbr
\eeq   Then, for each fixed $t$, $\Pt_{t}$ is a probability measure on $(\Omega,\sFt)$ and $\XtBt$ has a mean-zero BTBM density at $t$ on $\OFtPt$$:$ $\Pt_{t}[\XtBt\in dz]=\KBtzo dz$.
\eprop
\bpf 
Using the definition of $\Pt_{t}$ in \eqref{eq:RNexpl}, we have
\beq\lbl{eq:comdom}
\bsp
\Pt_{t}\lpa d\omega\rpa&=\XitmuXB\P\lpa d\omega\rpa=\exp\lbr\df{t}{\RBM}\lbk\mu\XB(t)-\frac12\mu^{2}t\rbk\rbr\P\lpa d\omega\rpa
\end{split}
\eeq
Thus, following Allouba and Zheng \cite{AZ01} and Allouba \cite{Abtp2}, we condition on the inner Brownian motion $B$ and we use the independence of $X$ and $B$ in the  BTBM $\XBt=X(|B_{t}|)$, to obtain
\beq\lbl{eq:comprange}
\bsp
\Pt_{t}\lbk \XBt\in dz\rbk
&=\lbk2\int_{0}^{\infty}e^{\frac{t}{s}\lpa\mu z-\frac12\mu^{2}t\rpa}\df{e^{-z^2/2s}}{\sqrt{2\pi s}}\df{e^{-s^2/2t}}{\sqrt{2\pi t}}ds\rbk dz
\\&=\lbk2\int_{0}^{\infty}\df{e^{-(z-\mu t)^2/2s}}{\sqrt{2\pi s}}\df{e^{-s^2/2t}}{\sqrt{2\pi t}}ds\rbk dz
\end{split}
\eeq
This means that
\ben\rencomrom
\item 
\beq\lbl{eq:Ptisprob}
\bsp
\Pt_{t}\lbk \Omega\rbk&=\int_{\R}\lbk2\int_{0}^{\infty}\df{e^{-(z-\mu t)^2/2s}}{\sqrt{2\pi s}}\df{e^{-s^2/2t}}{\sqrt{2\pi t}}ds\rbk dz=1,
\end{split}
\eeq
since the term in brackets is nothing more than the probability density $\KBtzmut$ of a BTBM $\XB^{\mu t}(t)$ with mean $\mu t$; and that
\item  under $\Pt_{t}$, $\XBt$ has the density of a BTBM at $t$ with mean $\mu t$; i.e., $$\Pt_{t}[\XBt\in dz]=\KBtzmut dz$$ 
(equivalently that $\XtBt=\XBt-\mu t$ is a zero-mean BTBM at $t$ with $\Pt_{t}[\XtBt\in dz]=\KBtzo dz$).  
\een
The proof is complete.
\epf

A quick computation, together with Fubini's theorem, gives us the fourth moment of a BTBM
\beq\lbl{eq:4thmo}
\bsp
\EP\X^{4}_{B}(t)&=\int_{\R}y^{4}\KBtyo dy\\
&=\int_{\R}y^{4}\lbk2\int_{0}^{\infty}\df{e^{-y^2/2s}}{\sqrt{2\pi s}}\df{e^{-s^2/2t}}{\sqrt{2\pi t}}ds\rbk dy
\\&=3\int_{0}^{\infty}\frac{ 2s^{2} {\mathrm e}^{-\frac{s^{2}}{2 t}}}{\sqrt{2\pi t}}ds=3t
\end{split}
\eeq
where the last integral is nothing more than the second moment of a zero mean normal random variable, with variance $t$.
This indicates that the quartic variation of $\XBt$ at $t$ is also $3t$ (see \prpref{prp:BTBMqv} and its proof below).  In addition, it's easy to obtain the following first $3$ moments of a mean-zero BTBM
\beq\lbl{eq:3mom}
\E|\XBt|=\frac{2^{\frac 5 4}\Gamma(\frac3 4)}{\pi}t^{1/4},\ \E|\XBt|^{3}=\frac{2^{\frac 5 4}}{\Gamma(\frac3 4)}t^{3/4}, \mbox{ and } \E|\XBt|^{2}=2\frac{\sqrt t}{\sqrt\pi}.
\eeq
\subsection{The BTBM quartic variation is $3t$: a new proof}
\bprop[The quartic variation of the BTBM is $3t$]\lbl{prp:BTBMqv}
Let $$\XB=\{\XBt,\sFt;t>0\}$$ be a BTBM on $\OFP$.  The
quartic variation of $\XB$ at $t$ is $\lang\XB\rang^{(4)}_t=3t$ almost surely.
\eprop
\brm\lbl{rem:Landallvariations}
\ben\rencomrom 
\item We note that, for the special case of the iterated Brownian motion (IBM), Burdzy proved the stronger result  of $L^{p}$ convergence---rather than convergence in probability---in the quartic variation expression \cite{BurdzyVar}.  On the other hand, our proof presented here is substantially shorter, and is easily adaptable to the whole class of Brownian-time processes in the Allouba and Zheng article \cite{AZ01}, which includes IBM as a special case. Moreover, with some adaptation of our proof, we also obtain $L^{2}$ convergence to the quartic variation $3t$.  This is largely due to our building on the fundamental conditioning strategy in Allouba and Zheng \cite{AZ01} and in Allouba \cite{Abtp2} in all of our proofs below.  
\item It follows from \lemref{lem:BTBMqv} together with  elementary stochastic analytic theory (e.g.Karatzas and Shreve \cite[p.~35]{KaratzasShreve}) that
\beq\lbl{eq:allvar}
\lang\XB\rang^{(p)}_t=\bc \infty, &\mbox{ for }p<4;\\
0,&\mbox{ for }p>4.\ec
\eeq
\een 
\erm

The proof of \prpref{prp:BTBMqv} proceeds via a couple of lemmas.  
\blm[Conditional expected fourth variation over $\Pi$ converges to $3t$]\lbl{lem:conpivar} Fix $t>0$ and let $\Pi$ and $|\Pi|$ be a partition of $[0,t]$ and its mesh, respectively, as in \defnref{def:pvar}.  Let $V_{t}^{(4)}(\XB,\Pi)$ be the fourth variation of the BTBM $\XB$ over $\Pi$ 
\beq\lbl{eq:4thvoverPi}
V_{t}^{(4)}(\XB,\Pi):=\sum_{k=1}^{n}\lab X(|B_{t_{k}}|)-X(|B_{t_{k-1}}|)\rab^{4}=\sum_{k=1}^{n}\lab\Delta\X_{k}\rab^{4}.
\eeq
If $\sGtB=\sigma\lpa B_s;0\le s\le t\rpa$ then \beq\lbl{eq:cond4thvar}
\E\lbk \left.V_{t}^{(4)}(\XB,\Pi)\rab\sG_{t}\rbk\xrightarrow [|\Pi|\searrow0]{\P}3t.
\eeq
\elm
\bpf
First, for any pair of finite collections of positive numbers $\{t_{1}, \ldots,t_{m}\}$ and $\{s_{1}, \ldots,s_{m}\}$, we have that the conditional finite dimensional distribution of the BTBM $\XB$ given the inner Brownian motion $B$ is
\beq\lbl{eq:findimdist}
\bsp
&\P\lbk\left.\XB(t_{1})\in dz_{1}, \ldots,\XB(t_{m})\in dz_{m}\right||B_{t_{1}}|=s_{1},\ldots,|B_{t_{m}}|=s_{m}\rbk
\\= \ &\P\lbk\left.X(|B_{t_{1}}|)\in dz_{1}, \ldots,X(|B_{t_{m}}|)\in dz_{m}\right||B(t_{1})|=s_{1},\ldots,|B(t_{m})|=s_{m}\rbk
\\= \ & \P\lbk\left.X(s_{1})\in dz_{1}, \ldots,X(s_{m})\in dz_{m}\right||B(t_{1})|=s_{1},\ldots,|B(t_{m})|=s_{m}\rbk
\\= \ & \P\lbk X(s_{1})\in dz_{1}, \ldots,X(s_{m})\in dz_{m}\rbk,
\end{split}
\eeq
where the last equality in \eqref{eq:findimdist} follows from the independence of the Brownian motions $X$ and $B$.  
This means that the conditional finite dimensional distributions of $\XB$ given the values of $B$ at the times $\{t_{1}, \ldots,t_{m}\}$ is the same as the finite dimensional distribution of the Brownian motion $X$.  By continuity of $X$ and $B$ this implies that the conditional law of $\XB$ given $B$ is the Gaussian law of $X$ with covariance 
\beq\lbl{eq:covX}
\bsp
&\Cov\lpa\left.\XB(t_{i}),\XB(t_{j})\right|\lab B(t_{i})\rab =s_{i},\lab B(t_{j})\rab =s_{j}\rpa
\\ &=\Cov\lpa X(s_{i}),X(s_{j})\rpa=s_{i}\wedge s_{j}.
\end{split}
\eeq

Equivalently, by definition of $\sGtB$ and from \eqref{eq:4thvoverPi}, we have that $B_{t_{k}} \in\sGtB$ for every $k=1,\ldots,n$. The independence of $B$ and $X$, and hence of $\sGtB$ and the outer Brownian motion $X$, implies that, the conditional law of the BTBM $\XB$, given the $\sigma$-field $\sGtB$, is centered Gaussian with covariance 
$$\Cov\lpa\left. X(|B_{u}|,X(|B_{v}|)\right|\sGtB\rpa=|B_{u}|\wedge |B_{v}|, \forall 0<u,v<\infty.$$  

Now, from \eqref{eq:4thvoverPi} we get
\beq\lbl{eq:scond4thvaroverpi}
\E\lbk \left.V_{t}^{(4)}(\XB,\Pi)\rab\sGtB\rbk=
\sum_{k=1}^{n}\E\lbk\left.\lab\Delta\X_{k}\rab^{4}\rab\sGtB\rbk.
\eeq
From the discussion above, given $\sGtB$, the increments $\lpa\Delta\X_{k}\rpa_{k=1}^{n}$ are centered jointly Gaussian.  Thus, using the standard fact that for the standard normal random variable $Z\eqdis N(0,\sigma^{2})$ we have $\E Z^{4}=\sigma^{4}$,  together with the trivial equality  $\Var(U-V) = \Var(U) + \Var(V) - 2\Cov(U,V)$, we have
\beq\lbl{eq:cond4thvaroverpi}
\bsp
\E\lbk\left.\lab\Delta\X_{k}\rab^{4}\rab\sG_{t}\rbk&=
3\lpa\E\lbk\left.\lab\Delta\X_{k}\rab^{2}\rab\sG_{t}\rbk\rpa^{2}\\
&=3\lpa|B_{t_{k}}|+|B_{t_{k-1}}|-2|B_{t_{k}}|\wedge|B_{t_{k-1}}|\rpa^{2}\\&= 3\lab |B_{t_{k}}| - |B_{t_{k-1}}|\rab^{2}.
\end{split}
\eeq
Thus, by \eqref{eq:scond4thvaroverpi} and \eqref{eq:cond4thvaroverpi} we get
\beq\lbl{eq:quadRBM}
\bsp
\lbk \left.V_{t}^{(4)}(\XB,\Pi)\rab\sG_{t}\rbk
&=\sum_{k=1}^{n}3\lab |B_{t_{k}}| - |B_{t_{k-1}}|\rab^{2}.
\\&\xrightarrow[|\Pi|\searrow0]{\P}3\lang |B| \rang^{(2)}_{t}
=3t,
\end{split}
\eeq
where the last equality follows from Tanaka's formula
\beq\lbl{eq:Tanaka}
|B_t|=\int_0^t \sgn(B_s)dB_s + L_t^0(B),
\eeq
where $L_t^0(B)$ is the local time of $B$ at $0$ (a continuous finite-variation process).  To see this, apply the quadratic variation to \eqref{eq:Tanaka} and use the properties of It\^o's stochastic integral and the fact that the quadratic variation of a finite variation process is zero to get 
\beqs
\bsp
\lang |B| \rang^{(2)}_{t}&=\lang \int_0^\cdot \sgn(B_s)dB_s + L_\cdot^0(B)\rang^{(2)}_{t}
\\&= \int_0^t \sgn^{2}(B_s)ds
=t.
\end{split}
\eeqs
We have proven \lemref{lem:conpivar}.
\epf
We now turn to the final major lemma.
\blm[Fourth variation over $\Pi$ converges to $3t$]\lbl{lem:BTBMqv}
Let $V_{t}^{(4)}(\XB,\Pi)$, $\Pi$, and $\sGtB$ be as in \lemref{lem:conpivar}. Suppose further that \eqref{eq:cond4thvar} in \lemref{lem:conpivar} holds.  Then,
\ig
\beq\lbl{eq:fvandconfv}
V_{t}^{(4)}(\XB,\Pi)-\E\lbk\left.V_{t}^{(4)}(\XB,\Pi)\right|\sGtB\rbk\xrightarrow[|\Pi|\searrow0]{\P}0
\eeq
and hence\fi
\beq\lbl{lem:QuartV} 
V_{t}^{(4)}(\XB,\Pi)\xrightarrow[|\Pi|\searrow0]{\P} 3t; 
\eeq
i.e., $\lang\XB\rang^{(4)}_t=3t$.
\elm
\bpf
The desired convergence \eqref{lem:QuartV} follows from \lemref{lem:conpivar} once we establish that
\beq\lbl{eq:condvarvan}
\Var\lpa\lftdt V_{t}^{(4)}(\XB,\Pi)\rab \sGtB\rpa\xrightarrow[|\Pi|\searrow0]{\P}0.
\eeq
This is easily seen by the following argument.  We assume that \eqref{eq:condvarvan} holds.  Clearly, we have 
\beq\lbl{eq:2parts}
\bsp
&\P\lpa\lab V_{t}^{(4)}(\XB,\Pi) - 3t\rab > \epsilon\rpa\\ &\le \P\lpa\lab V_{t}^{(4)}(\XB,\Pi) - \E\lbk\lftdt V_{t}^{(4)}(\XB,\Pi)\rab\sGtB\rbk\rab > \frac{\epsilon}{2}\rpa 
\\&
+ \P\lpa\lab\E\lbk\lftdt V_{t}^{(4)}(\XB,\Pi)\rab \sGtB\rbk - 3t\rab > \frac{\epsilon}{2}\rpa.
\end{split}
\eeq
To deal with the first term on the right side of the inequality of \eqref{eq:2parts}, we condition on $\sGtB$ then take the expectation.  Doing this gives us
\beq\lbl{eq:1sttermrhs0}
\bsp
&\P\lpa\lab V_{t}^{(4)}(\XB,\Pi) - \E\lbk\lftdt V_{t}^{(4)}(\XB,\Pi)\rab\sGtB\rbk\rab > \frac{\epsilon}{2}\rpa
\\&= \E\left[ \P\lftdt\lpa\lab V_{t}^{(4)}(\XB,\Pi) - \E\lbk\lftdt V_{t}^{(4)}(\XB,\Pi)\rab\sGtB\rbk \rab > \frac{\epsilon}{2} \rab \sGtB\right)\right].
\end{split}
\eeq
Applying the conditional Chebyshev inequality to \eqref{eq:1sttermrhs0} we obtain
\beq\lbl{eq:cndCheb}
\bsp
\P&\lftdt\lpa\lab V_{t}^{(4)}(\XB,\Pi) - \E\lbk\lftdt V_{t}^{(4)}(\XB,\Pi)\rab\sGtB\rbk \rab > \frac{\epsilon}{2} \rab \sGtB \rpa
\\&\leq \frac{4\Var\lpa\lftdt V_{t}^{(4)}(\XB,\Pi)\rab\sGtB\rpa}{\epsilon^2}\wedge1\xrightarrow[|\Pi|\searrow0]{\P}0,
\end{split}
\eeq
where the convergence to zero is from our assumption in \eqref{eq:condvarvan}.  Taking the expected value on both sides of \eqref{eq:cndCheb}, along with \eqref{eq:1sttermrhs0}, and the dominated convergence theorem we get
\beq\lbl{eq:1sttermrhs}
\bsp
&\P\lpa\lab V_{t}^{(4)}(\XB,\Pi) - \E\lbk\lftdt V_{t}^{(4)}(\XB,\Pi)\rab\sGtB\rbk\rab > \frac{\epsilon}{2}\rpa
\\&\le\E\lbk\frac{4\Var\lpa\lftdt V_{t}^{(4)}(\XB,\Pi)\rab\sGtB\rpa}{\epsilon^2}\wedge1\rbk\xrightarrow[|\Pi|\searrow0]{}0. 
\end{split}
\eeq
As for the second term on the right side of the inequality in \eqref{eq:2parts}, we have by assumption that \eqref{eq:cond4thvar} holds; i.e.,
\beq\lbl{eq:cond4thvar3t}
\P\lpa\lab\E\lbk \left.V_{t}^{(4)}(\XB,\Pi)\rab\sG_{t}\rbk-3t\rab>\frac{\epsilon}{2}\rpa\xrightarrow [|\Pi|\searrow0]{}0.
\eeq
Our claim now follows by applying both \eqref{eq:1sttermrhs} and \eqref{eq:cond4thvar3t} to \eqref{eq:2parts}.

Lastly, it remains for us to prove that the convergence to zero in \eqref{eq:condvarvan} holds. Let's start by observing that the conditional variance
$$\Var\lpa\lftdt V_{t}^{(4)}(\XB,\Pi)\rab \sGtB\rpa$$ has the following variance-covariance decomposition 
\beq\lbl{eq:condvar}
\bsp
\Var\lpa\lftdt V_{t}^{(4)}(\XB,\Pi)\rab \sGtB\rpa
&= \sum_{k=1}^{n}\Var\bigl((\Delta \mathbb{X}_k)^4\mid \sGtB\bigr)\\
  &+ \sum_{i\neq j}\Cov\bigl((\Delta \mathbb{X}_i)^4,(\Delta \mathbb{X}_j)^4\mid \sGtB\bigr).
\end{split}
\eeq
Now, given $\sGtB$, $\Delta \mathbb{X}_k\eqdis N(0,\sigma_k^2)$ with $
\sigma_k^2=\bigl||B_{t_{k}}|-|B_{t_{k-1}}|\bigr|$ (the size of the $k$-th Brownian-time increment).   Thus, the variance (diagonal) term in \eqref{eq:condvar} satisfies the estimate
\beq\lbl{eq:diagsumvan}
\bsp
\sum_{k=1}^{n}\Var\bigl((\Delta \mathbb{X}_k)^4\mid \sGtB\bigr)
&=96\sum_{k=1}^{n}\sigma_k^8\\
&=96\sum_{k=1}^{n}\bigl||B_{t_{k}}|-|B_{t_{k-1}}|\bigr|^4
\\&\le\max_{1\le k\le n}\bigl||B_{t_{k}}|-|B_{t_{k-1}}|\bigr|^2
\sum_{k=1}^{n}\bigl||B_{t_{k}}|-|B_{t_{k-1}}|\bigr|^2
\\ &\xrightarrow[|\Pi|\searrow0]{\P}0,
\end{split}
\eeq
where we have used the fact that if $Z\eqdis N(0,\sigma^2)$, then $\mathbb{E}[Z^8]=105\sigma^8$.
Hence $\Var(Z^4)=105\sigma^8-(3\sigma^4)^2=96\sigma^8$.
The last convergence to zero assertion follows by (i) the uniform continuity of $B$ on $[0,t]$, hence 
$$\max_{1\le k\le n}\bigl||B_{t_{k}}|-|B_{t_{k-1}}|\bigr|^2
\xrightarrow[|\Pi|\searrow0]{\mbox{a.s.}}0,$$  
and (ii) by the fact that 
\beqs
\bsp
\sum_{k=1}^{n}\bigl||B_{t_{k}}|-|B_{t_{k-1}}|\bigr|^2 &\xrightarrow[|\Pi|\searrow0]{\P}\lang|B|\rang_{t}^{(2)}=t.\end{split}
\eeqs
We have proven that the conditional variance term in \eqrf{eq:condvar} vanishes in probability.

We now turn to the covariance (off-diagonal) term of \eqrf{eq:condvar}.  
First, for each pair $i\ne j,$ we apply the  Cauchy-Schwarz inequality for conditional covariance, to get  
\beq\lbl{eq:condCS}
\bsp
&\lab\Cov\bigl((\Delta \mathbb{X}_i)^4,(\Delta \mathbb{X}_j)^4\mid \sGtB\bigr)\rab
\\&\le \sqrt{\Var\bigl((\Delta \mathbb{X}_i)^4\mid \sGtB\bigr)}\sqrt{\Var\bigl((\Delta \mathbb{X}_j)^4\mid \sGtB\bigr)},
\end{split}
\eeq
Again, given $\sGtB$, the increment $\Delta \X_k = X(|B_{t_k}|) - X(|B_{t_{k-1}}|)$ is normally distributed with variance $\sigma^{2}_k = ||B_{t_k}| - |B_{t_{k-1}}||:=\delta_{k}$.  As we observed above $\Var(Z^{4})=96\sigma_{k}^{8}=96\delta_{k}^{4}$ for any $Z\eqdis N(0,\sigma_{k}^{2})$.  So, in terms of the Brownian-time increments $\delta_{k}$, the Cauchy-Shwarz bound in \eqref{eq:condCS} becomes 
\beq\lbl{eq:condCS2}
\bsp
\lab\Cov\bigl((\Delta \mathbb{X}_i)^4,(\Delta \mathbb{X}_j)^4\mid \sGtB\bigr)\rab
\le 96 \delta_i^2 \delta_j^2.
\end{split}
\eeq
We now observe that the covariance is exactly zero if the Brownian-time subintervals 
\beqs
\bsp
I_i &:= \lbk|B_{t_{i-1}}|\wedge |B_{t_i}|, |B_{t_{i-1}}|\vee |B_{t_i}|)\rbk \mbox{ and }\\ 
I_{j}&:= \lbk|B_{t_{j-1}}|\wedge |B_{t_j}|, |B_{t_{j-1}}|\vee |B_{t_j}|)\rbk
\end{split}
\eeqs
are disjoint.  Thus,  
\beq\lbl{eq:condCS3}
\bsp
\sum_{i\ne j}\lab\Cov\bigl((\Delta \mathbb{X}_i)^4,(\Delta \mathbb{X}_j)^4\mid \sGtB\bigr)\rab
\le 96 \sum_{i\ne j}\delta_i^2 \delta_j^2\ind_{\{I_i \cap I_j \neq \emptyset\}}.
\end{split}
\eeq   For small $|\Pi|$ the Brownian-time increments $\delta_{k}$ are small, and 
\beq\lbl{eq:smallincr}
\bsp
&\sum_{i\ne j}\lab\Cov\bigl((\Delta \mathbb{X}_i)^4,(\Delta \mathbb{X}_j)^4\mid \sGtB\bigr)\rab
\\&\le
96\lpa\max_{1\le k\le n}\delta_{k}^{2}\rpa\sum_{i,j=1}^{n} \delta_i \delta_j \cdot \ind_{\{I_i \cap I_j \neq \emptyset\}}
\end{split}
\eeq
Intuitively, the sum on the right side of \eqref{eq:smallincr} counts how often two increments overlap in level space.  In the limit, as we shall see below, this sum becomes the integral of square of the local time (the self intersection measure). 
To deal effectively with the sum on the right side of \eqref{eq:smallincr}, we briefly recall the definition of local time for the reflecting Brownian motion $|B|$ at the point $a$, ${L^{a}_{t}} $:
\beq\lbl{eq:locT}
\bsp
{L^{a}_{t}}&=\lim_{\epsilon\searrow0}\frac1{2\epsilon}\int_{0}^{t}\ind_{\{a-\epsilon<|B_{s}|<a+\epsilon\}}ds, 
\end{split}
\eeq
where the limit above always exist in the almost sure sense, The basic idea is that 
$L^{a}_{t} $ is an (appropriately rescaled and time-parametrized) measure of how much time 
$|B_{s}|$ has spent at $a$ up to time $t$.
The occupation time of region $\mathcal R\in\sB(\Rp)$, by $\lbr|B_{s}|,0\le s\le t\rbr$, is thus given, in terms of the occupation measure $\mu_{t}$, by
\beq\lbl{eq:OccT}
\mu_{t}(\mathcal R) = \int_0^t \ind_{\{|B_s| \in \mathcal R\}} ds = \int_{\mathcal R} L_t^a da;
\eeq
i.e., the local time $L_t^a$ is the density of $\mu_{t}$ with respect to Lebesgue measure: $\mu_{t}(da)=L_t^a da$.  

Observe that
\[
\delta_i \delta_j \mathbf{1}_{\{I_i \cap I_j \neq \emptyset\}}
= \int_0^\infty \mathbf{1}_{\{a \in I_i\}} \mathbf{1}_{\{a \in I_j\}} \, da.
\]
Hence,
\beq\lbl{eq:intersum}
\sum_{i,j=1}^n \delta_i \delta_j \mathbf{1}_{\{I_i \cap I_j \neq \emptyset\}}
= \int_0^\infty \left( \sum_{i=1}^n \mathbf{1}_{\{a \in I_i\}} \delta_i \right)^2 \, da.
\eeq
By Tanaka's formula, the local time $L_t^a$ of the reflected BM $|B|$ at level $a$ satisfies
\[
L_t^a = \lim_{|\Pi| \to 0} \sum_{i=1}^n \mathbf{1}_{\{a \in I_i\}} \delta_i,
\]
where this convergence holds almost surely (a.s.) for each fixed $a$ because the sum is a pathwise Riemann-Stieltjes approximation of the local time.  
We now have
\beq\lbl{eq:limitofsum}
\bsp
\lim_{|\Pi| \to 0}\sum_{i,j=1}^n \delta_i \delta_j \mathbf{1}_{\{I_i \cap I_j \neq \emptyset\}}
&= \lim_{|\Pi| \to 0}\int_0^\infty \left( \sum_i \mathbf{1}_{\{a \in I_i\}} \delta_i \right)^2 da
\\&=\int_0^\infty \lim_{|\Pi| \to 0}\left( \sum_i \mathbf{1}_{\{a \in I_i\}} \delta_i \right)^2 da
\\&=\int_0^\infty (L_t^a)^2 da.
\end{split}
\eeq
Since $a \mapsto L_t^a$ is continuous with compact support (because Brownian motion has compact range on $[0,t]$), we can pass the limit inside the integral using the dominated convergence theorem. Therefore\footnote{The convergence is almost sure. It is also true in $L^p$ for $1 \le p < \infty$ due to bounded moments of local time, but the natural formulation is almost sure convergence. },

\[
\sum_{i,j=1}^n \delta_i \delta_j \mathbf{1}_{\{I_i \cap I_j \neq \emptyset\}}
\;\xrightarrow[|\Pi|\to 0]{\text{a.s.}}\;
\int_0^\infty (L_t^a)^2  da.
\]

From \eqref{eq:smallincr} and \eqref{eq:convtoselfinters} we have, for small mesh size $|\Pi|$, that
\beq\lbl{eq:convto0}
\bsp
&\sum_{i\ne j}\lab\Cov\bigl((\Delta \mathbb{X}_i)^4,(\Delta \mathbb{X}_j)^4\mid \sGtB\bigr)\rab
\\&\le C \lpa\max_{1\le k\le n}\delta_{k}\rpa\int_{\Rp}\lpa L_{t}^{a}\rpa^{2} da
\xrightarrow[|\Pi|\searrow0]{\mbox{a.s.}}0,
\end{split}
\eeq
where the last convergence-to-zero assertion follows from the almost sure continuity of $|B|$, so $\lim_{|\Pi|\searrow0} \max_{1\le k\le n}\delta_{k}=0$ together with the fact that, for each $t$, the local time $L_{t}^{a}$ is continuous and has compact support in the spatial variable $a$.  Thus,
\beq\lbl{eq:convto0P}
\bsp
\sum_{i\ne j}\lab\Cov\bigl((\Delta \mathbb{X}_i)^4,(\Delta \mathbb{X}_j)^4\mid \sGtB\bigr)\rab
\xrightarrow[|\Pi|\searrow0]{\P}0,
\end{split}
\eeq
as desired.
\epf
\bpfs{Proof of \prpref{prp:BTBMqv}} The proof is now complete by \lemref{lem:conpivar} and \lemref{lem:BTBMqv} above.
\epfs
We have reached the last part of the proof of our main result.
\bpfs{Proof of \thmref{thm:BTCOM}} The proof is now complete by \prpref{lem:crux} and \prpref{prp:BTBMqv}.
\epfs

\
\section{Appendix: Convergence of the overlap sum to the self intersection meassure via delta functions}
Here, we present another intuitively appealing delta function argument for the convergence of the overlap sum  $\sum_{i,j} \delta_i \delta_j \ind_{\{I_i \cap I_j \neq \emptyset\}}$ to the intgral of the square of the local time $L_{t}^{a}$ of the reflected BM $|B|$.  First, recall the Dirac delta function definition of local time for the reflected Brownian motion $|B|$ at the point $a$, ${L^{a}_{t}} $:
\beq\lbl{eq:locTdelta}
\bsp
{L^{a}_{t}}&= \int_{0}^{t}{\mathfrak d}_{a}(|B_{s}|)ds, 
\end{split}
\eeq
and where we used ${\mathfrak d}_{a}(|B_{s}|)$ to denote the Dirac measure of the path of $|B_{s}|$, concentrated at point $a$\footnote{Of course, we may think of  ${\mathfrak d}_{a}(|B_{s}|)$ as Dirac delta function at $|B_{s}|$, centered at $a$.}.

\beq\lbl{eq:inters1}
\bsp
\sum_{i,j} \delta_i \delta_j \ind_{\{I_i \cap I_j \neq \emptyset\}}\xrightarrow[|\Pi|\searrow0]{\mbox{a.s.}}\int_{0} ^{t} \int_{0} ^{t}\dirz(|B_s| - |B_r|) ds dr,
\end{split}
\eeq
where the right side term is the integral of the measure of the difference of the paths.  Integrating over all possible intersection points $a\in\Rp$ and using the local time definition \eqref{eq:locTdelta}, we have
\beq\lbl{eq:dintprod}
\bsp
\int_{0} ^{t} \int_{0} ^{t}\dirz(|B_s| - |B_r|) ds dr
&=\int_{\Rp}\lpa\int_{0} ^{t}\dira(|B_s|) ds\rpa\lpa\int_{0} ^{t}\dira(|B_r|) dr\rpa da
\\&=\int_{\Rp}\lpa L_{t}^{a}\rpa^{2} da.
\end{split}
\eeq
From \eqref{eq:inters1} and \eqref{eq:dintprod} we obtain the well known self-intersection measure for the reflected BM $|B|$ as the limit:
\beq\lbl{eq:convtoselfinters}
\sum_{i,j} \delta_i \delta_j \ind_{\{I_i \cap I_j \neq \emptyset\}}\xrightarrow[|\Pi|\searrow0]{\mbox{a.s.}}\int_{\Rp}\lpa L_{t}^{a}\rpa^{2} da.
\eeq
\section{Frequent acronyms and notations key}\lbl{sec:acrnot}
\begin{enumerate}\renewcommand{\labelenumi}{\Roman{enumi}.}
\item {\textbf{Acronyms}}\vspace{2mm}
\begin{enumerate}\renewcommand{\labelenumii}{(\arabic{enumii})}
\item BM: Brownian motion.
\item BTBM: Brownian-time Brownian motion.
\end{enumerate}
\vspace{2.5mm}
\item {\textbf{Notations}}\vspace{2mm}\\
\ig
Throughout the article sets of numbers or Euclidean spaces are denoted using the {\tt{$\backslash$mathbb}} font  (e.g., $\R$, $\Rd$, $\N$, etc...). \fi The following are some of the notations we use in this article:
\begin{enumerate}\renewcommand{\labelenumii}{(\arabic{enumii})}
\item $\KBtxy$ is the density or kernel of the BTBM, starting at $x$ (with mean x);
\ig
\item $\Wtt(\Rd;\R)$ is the Sobolev space of $\Ltwo(\Rd;\R)$ functions whose weak derivative up to order $2$ are also $\Ltwo(\Rd;\R)$;
\item $\W^{-2,2}(\Rd;\R)$ is the dual of $\Wtt(\Rd;\R)$; 
\item $\pa_{\l}$ is the partial derivative in the $\l$-th spatial variable, $\l\in\{1,\ldots,d\}$;
\item $\pa^{2}_{\l^{2}}$ is the second partial derivative in the $\l$-th spatial variable, $\l\in\{1,\ldots,d\}$;
\item $\pa_{t}$ is the partial derivative in the time variable $t$;
\fi
\item $\KBtzo$ denotes the BTBM density starting at $0$ (with zero mean);
\end{enumerate}
\end{enumerate}

\end{document}